\documentclass {amsart}
\usepackage{amsmath, amssymb}
\usepackage{amssymb,pb-diagram}
 \usepackage{amscd}

\newcommand{\ZZ}{\mathbb Z}
\newcommand{\PP}{\mathbb P}

\newcommand{\CC}{\mathbb C}

\newcommand{\caC}{{\mathcal C}}
\newcommand{\caX}{{\mathcal X}}

\newcommand{\Aut}{\mathop {\rm {Aut}}\nolimits}

\newtheorem{thm}{Theorem}[section]
\newtheorem{cor}{Corollary}[section]

\newtheorem{lem}{Lemma}[section]

\newtheorem{exmple}{Example}[section]
\newtheorem{rem}{Remark}[section]
\newtheorem{qz}{Question}[section]
\newcommand{\proofend}{\qed \par\smallskip\noindent}

\renewcommand{\thesubparagraph}{\theparagraph.\@arabic\c@subparagraph}
\usepackage{ifpdf}
\ifpdf
	\usepackage[pdftex]{graphicx}
	\DeclareGraphicsExtensions{.pdf}
	\usepackage{hyperref}
\else
	\usepackage[dvips]{graphicx}
	\DeclareGraphicsExtensions{.eps}
\fi

\begin{document}

\title{A note on local trigonal fibrations
 }

\dedicatory{To the memory of  Professor Ngyuen Huu Duc.}
\author{Mizuho ISHIZAKA and Hiro-o TOKUNAGA}\footnote{MSC 2000 Classfication: primary 14H15, secondary 14H30, 14H51, Key words: 
trigonal curve, triple cover}
\maketitle

%

%
\large
\textbf{Introduction}\label{intro}
\normalsize
\par\bigskip

   Let $\varphi : \caC =\{\caC_t\}_{t\in \Delta_{\epsilon}} \to \Delta_{\epsilon}$ be a local family of curves over a
   small disc $\Delta_{\epsilon} := \{t \in \CC\mid |t| < \epsilon\}$, 
   $\epsilon$: a small positive real number,  which means that $\caC$ is a smooth
   surface, $\varphi$ is proper and surjective, all fibers of $\varphi$ is
   connected and $\varphi^{-1}(t)$ is smooth for $t\neq 0$.    We call $\varphi : \caC \to \Delta_{\epsilon}$ a local hyperelliptic (resp. trigonal) fibration of genus $g$ if
   $\caC_t$ $(t \neq 0)$ is a hyperelliptic (resp. trigonal) curve of genus $g$.
   
   In \cite{chen-tan}, Chen and Tan gave two examples of local trigonal fibrations of genus $3$ such that
   their central fibers, $\caC_0$, are smooth hyperelliptic curves of genus $3$.
     A purpose of this note is to show that any hyperelliptic curve of genus $g$ appears as the
   central fiber of a local trigonal fibration.

  Let us explain our setting and question precisely.   
  Let $B_1$ be a reduced  divisor on $\PP^1\times \Delta_{\epsilon}$
  such 
  that $B_1$ meets $\PP^1\times \{t\} (t\neq 0)$ at $2g+2$ distinct points transversely.  In particular, all singularities of $B_1$ are in $\PP^1\times\{0\}$.
   Let $f_o : W \to \PP^1\times\Delta_{\epsilon}$ be the double
  cover with $\Delta_{f_o} = B_1$, where $\Delta_{f_o}$ denotes the branch locus of $f_o$,
  and let $\mu_o : \widetilde{W} \to W$ be the canonical resolution (see
  \cite{horikawa} for the canonical resolution).  
 We consider a local hyperelliptic fibration of genus $g$ given by putting $\caC = \widetilde{W}$ and
  $\varphi_o:= \mbox{\rm pr}_2\circ f_o \circ \mu_o$. 
  
   We say that 
   
   \medskip
   
   $(i)$ $\varphi_o : \caC \to \Delta_{\epsilon}$ is \textit{of horizontal type}, if $B_1$ does not contain $\PP^1\times \{0\}$ as its
   irreducible component,  and
   
   \medskip
   
   $(ii)$ $\varphi_o : \caC \to \Delta_{\epsilon}$ is \textit{of non-horizontal type}, if 
   $B_1$ contains $\PP^1\times \{0\}$ as its
   irreducible component.
   
   \medskip
   
   Note that we may assume that $\varphi_o$ is either horizontal or non-horizontal by taking
   $\epsilon (> 0)$ small enough. Also any local hyperelliptic fibration is obtaned as the 
   relatively minimal model of $\varphi_o : \caC \to \Delta_{\epsilon}$. We now formulate
   our question as follows:

   \begin{qz}\label{qz:mondai} {\rm For an arbitrary local hyperelliptic fibration $\varphi_o : \caC
   \to \Delta_{\epsilon}$, does there exist a local trigonal fibration $\varphi : \widetilde{\caC}
   \to \Delta_{\epsilon}$ such that $\caC_0 = \widetilde{\caC}_0$?
   In other words, can the central fiber of any local hyperelliptic fibration 
   appear as that of a certain local trigonal fibration?
   }
   \end{qz}

  In this note, we give an answer to Question~\ref{qz:mondai} in the case when
  $\varphi_o$ is of horizontal type.
 
  \begin{thm}\label{thm:main}{Question~\ref{qz:mondai} is true for a local hyperelliptic
  fibration of horizontal type.
  }
  \end{thm}
  
  Since any hyperelliptic curve of genus $g$ appears as the centeral fiber of a local
  hyperelliptic fibration of horizontal type, i.e., a trivial fibration,  we have
  
   \begin{cor}\label{cor:co1r}{Let $D$ be any hyperelliptic curve of genus $g$. There exists a local trigonal
   fibration $\widetilde{\varphi}: \widetilde{\caC} \to \Delta_{\epsilon}$ such that $\widetilde{\caC}_0 = D$. 
   }
   \end{cor}
   
%
%

%
%

\section{Covers and the fundamental group}\label{sec:junbi}
Let $X$ and $M$ be a normal variety and a complex manifold, respectively. We call $X$ a \emph{(branched) cover of $M$} 
if there exists a finite surjective morphism $\pi : X \to M$. 
When needed, the covering morphism will be specified as \emph{a cover $\pi : X \to Y$}.
Let $B$ be a reduced divisor on $M$. The following facts are well-known:

\begin{itemize}

\item Choose a point $\ast \in M\setminus B$. Then the inclusion morphism $\iota : M\setminus B
\to M$ induces an epimorphism $\iota_* : \pi_1(M\setminus B, \ast) \to \pi_1(M, \ast)$.

\item Let $H$ be a subgroup of $\pi_1(M\setminus B, \ast)$. Then there exists a unramified cover $\pi_H : X_o \to M\setminus B$ over $M\setminus B$ with
$\pi_H(X_o, \hat{\ast}) \cong H$, $\pi_H(\hat{\ast}) = \ast$, such that
$X_o$ can be extened over $M$ uniquely. 
We also denote the extended cover of $\pi_H : X_0 \to M\setminus B$ 
by $\pi_H : X \to M$.
Note that the branch locus $\Delta_{\pi_H}$ of $\pi_H$ is a subset of $B$ and
$\deg \pi_H = [\pi_1(M\setminus B, \ast): H]$.
Conversely, if there exists a branched cover $\pi : X \to M$ with $\deg \pi = n$, then
there exists  a subgroup $H_{\pi}$ of $\pi_1(M\setminus \Delta_{\pi}, \ast)$ of 
index  $n$.

\item Let $H$ be a normal subgroup of $\pi_1(M\setminus B, \ast)$. Then there exists a 
Galois cover $\pi_H :  X \to M$ with $\Aut_M(X) \cong \pi_1(M\setminus B, \ast)/H$.

\end{itemize}

For the first statement, see \cite{shimada-tokunaga}, for example and for the last two statements, see \cite{sga} EXPOSE XII, for example.

Let $G$ be a finite group. We call a Galois cover with $\Aut_M(X)\cong G$ a \textit{$G$-cover}.
\medskip

The following lemma is fundamental throughout this article.

\begin{lem}\label{lem:fund}{Let $M$ and $B$ as above. Let
$D_6$ be the dihedral group of order $6$, which we describe  by $\langle \sigma, \tau \mid \sigma^2=\tau^3 = (\sigma\tau)^2 = 1\rangle$.
If there exists a non-Galois triple cover $\pi : X \to M$ with
$\Delta_{\pi}=B$, then there exists a $D_6$-cover $\hat {\pi} : \hat{X} \to M$,
with $\Delta_{\hat{\pi}} = B$. Conversely, if 
there exists a $D_6$-cover $\hat {\pi} : \widehat{X} \to M$
with $\Delta_{\hat{\pi}}=B$, then the quotient surface
$X := \widehat{X}/\langle\sigma\rangle$ by $\sigma$ admits a non-Galois triple cover $\pi : X \to M$ with
$\Delta_{\pi} = B$.
}
\end{lem}

\proof Let $\pi : X \to M$ be a given non-Galois triple cover of $M$ as above. Then there exists
a corresponding subgroup $H$ of $\pi_1(M\setminus B, \ast)$ with $\pi = \pi_H$ and
$[\pi_1(M\setminus B, \ast) : H] = 3$. Since $\pi$ is not Galois, $H$ is not normal. This 
implies that there exists an epimorphism $\pi_1(M\setminus B, \ast) \to D_6$ such that $K:=
\ker(\pi_1(M\setminus B, \ast) \to D_6) \subset H$. The $D_6$-cover corresponding to $K$ satisfies the desired property. Conversely, suppose that
there exists a $D_6$-cover $\varpi : \caX \to M$ with $\Delta_{\varpi} = B$.
As we have an epimorphism $\pi_1(M\setminus B, \ast) \to D_6$, there exists
a subgroup, $H$, of $\pi_1(M\setminus, B, \ast)$ such that $[\pi_1(M\setminus B, \ast) : H]= 3$ and $H$ contains $\ker(\pi_1(M\setminus B, \ast) \to D_6)$. The triple
cover corresponding to $H$ is the desired one. As for the equality
$\Delta_{\varpi} = \Delta_{\pi_H}$, see \cite{tokunaga2}. \proofend

\begin{rem}\label{cor:fund}{\rm Let $\pi : X \to M$ be as in Lemma \ref{lem:fund}.
There exists a double cover corresponding to the preimage of 
the subroup of order $3$ in $D_6$. We denote it by $\beta_1(\pi) : D(X/M) \to M$. Note that
$\Delta_{\beta_1(\pi)} \subseteq \Delta_{\pi}$.
}
\end{rem}


\section{Proof of Theorem\ref{thm:main}}

\subsection{Settings}\label{subsec:settings}

Throughout this section, we always assume that 

$(\ast)$ a local hyperelliptic fibration is of horizontal type.

Let $B_1$ be the branch locus of $f_o : W \to \PP^1\times \Delta_{\epsilon}$ as
in Introduction.
Let $B_o$ be a double section of 
$\mathrm{pr}_2 : \PP^1\times\Delta_{\epsilon} \to \Delta_{\epsilon}$
such that
\begin{itemize}

\item $B_o$ is smooth,

\item $\mathrm{pr}_2|_{B_o} : B_o \to \Delta_{\epsilon}$ has a unique ramification point over $0 \in \Delta_{\epsilon}$, and

\item $B_o \cap B_1 = \emptyset$.

\end{itemize}

Put $B := B_o + B_1$ and 
$G_B := \pi_1(\PP^1\times\Delta_{\epsilon}\setminus B, \ast)$. In order to study
non-Galois triple covers with branch locus $B$, 
by Lemma~\ref{lem:fund}, we need to know the description of
$G_B$ and its normal subgroup $K$ with $G_B/K \cong D_6$.

\subsection{A description of $G_B$}\label{subsec:description}

We describe $G_B$ via generators and their relation. The method used here
is well-known in computing the fundamental group of the complement of plane
curve via so-called ``Zariski-van Kampen method." We refer \cite{acc, act, dimca, 
shimada-tokunaga} and use results there freely.

Choose a point $t_o$ in $\Delta_{\epsilon}^*:= \Delta_{\epsilon}\setminus \{0\}$.
 Put $\PP^1_{t_o}:=\PP^1\times\{t_o\}$ and $\mathbf{Y}:=
\PP^1_{t_o}\cap B$. Let $\eta$ be a loop given by
$\{|t_o|\exp({2\pi i\theta})\mid 0 \le \theta \le 1\}$. Note that
$\pi_1(\Delta_{\epsilon}^*, t_o) \cong \langle \eta \rangle \cong \ZZ$.
Choose a point $b =(\hat{t}_o, t_o) \in (\PP^1\times\Delta_{\epsilon}^*)\setminus B$
with $(\{\hat{t}_o\}\times \Delta_{\epsilon}) \cap B = \emptyset$ 
and a geometric 
basis $\gamma_1, \gamma_2, \ldots, \gamma_{2g+4}$ of
$\pi_1(\PP^1_{t_o}\setminus \mathbf{Y}, \hat{t}_o)$ (see \cite[Definition 1.13]{acc} for
the definition of  a geometric basis). One can define a right action of
$\pi_1(\Delta_{\epsilon}^*, t_o)$ on 
$\pi_1(\PP^1_{t_o}\setminus \mathbf{Y}, \hat{t}_o)$
and $G_B$ is described through this action as follows:
\[
G_B \cong \langle \gamma_1, \ldots, \gamma_{2g+4}\mid \gamma_i^{\eta}
= \gamma_i \, \, (i =1,\ldots, 2g+4),\,\, \gamma_1\ldots\gamma_{2g+4} =1 \rangle.
\]
We may assume that

$\gamma_1, \ldots, \gamma_{2g+2}$ are meridians for the points in
$\PP_{t_o}^1\cap B_1$,

\noindent and

$\gamma_{2g+3}$ and $\gamma_{2g+4}$ are meridians for $\PP^1_{t_o}\cap
B_o$.

\medskip

Under these circumstances,  we have

\begin{lem}\label{lem:fundametal_group}{
\[
G_B \cong \langle \gamma_1,\ldots, \gamma_{2g+4} \mid \gamma_i^{\eta} =
\gamma_i \,\, (i = 1,\ldots, 2g+2), \,\, \gamma_{2g+3} = \gamma_{2g+4}, \,\,
\gamma_1\cdots\gamma_{2g+4} = 1 \rangle.
\]
}
\end{lem}

Let $H$ be the subgroup of $G_B$ generated by 
\begin{center}

$\gamma_i^2\,\, (i =1,\ldots, 2g+4)$ and $\gamma_j\gamma_{j+1}\,\,
(j = 1,\ldots, 2g+3)$.

\end{center}
Let $f: S \to \PP^1\times\Delta_{\epsilon}$ be  a double cover with $\Delta_f = B$.

\begin{lem}\label{lem:key-1}{$H$ is the subgroup of $G_B$ corresponding
to the double cover
\[
f':=f|_{S\setminus f^{-1}(B)} : S\setminus f^{-1}(B) \to 
(\PP^1\times \Delta_{\epsilon})\setminus B.
\]
}
\end{lem}

\proof Put $\overline{\varphi}:= \mathrm{pr}_2\circ f$ and 
$S_{t_o} := \overline{\varphi}^{-1}(t_o)$. 
Let $f_{t_o} : S_{t_o} 
\to \PP^1_{t_o}$  be the restriction of $f$ 
to the fiber over $t_o$ and put $f'_{t_o}:=f_{t_o}|_{S_{t_o} \setminus f_{t_o}^{-1}(\mathbf{Y})}$. 
We have a commutative diagram:
\[
\begin{CD}
S_{t_o}\setminus f_{t_o}^{-1}(\mathbf{Y}) @>>> S\setminus f^{-1}(B) \\
@V{f'_{t_o}}VV                 @VV{f'}V \\
\PP^1_{t_o}\setminus \mathbf{Y} @>>> (\PP^1\times\Delta_{\epsilon})\setminus B.
\end{CD}
\]
From the above diagram, we have a commutative diagram of groups:
\[
\begin{CD}
\widetilde{N} @>>> N \\
@V{(f'_{t_o})_{\sharp}}VV                 @VV{(f')_{\sharp}}V \\
\widetilde{G}:= \pi_1(\PP^1_{t_o}\setminus \mathbf{Y}, \hat{t}_o) @>>> G_B,
\end{CD}
\]
where $\widetilde{N}$ and $N$ are the subgroups of index $2$ corresponding
to  the doube covers $f'_{t_o}$ and $f'$. Note that $\widetilde{N} \to N$ is 
surjective by \cite[Theorem 2.30]{shimada-tokunaga} and $[G_B:N] = 2$. 
Our statement follows from the claim below:

\textbf{Claim} \textit{Let $\widetilde{H}$ be the subgroup of $\widetilde{G}$
 generated by 
\begin{center}
$\gamma_i^2\,\, (i =1,\ldots, 2g+4)$ and $\gamma_j\gamma_{j+1}\,\,
(j = 1,\ldots, 2g+3)$.
\end{center}
Then $\widetilde{H} = \widetilde{N}$.
}

\medskip

\noindent \textsl{Proof of Claim.} Since the branch locus of $f_{t_o} :
S_{t_o} \to \PP^1_{t_o}$ is $\mathbf{Y}$,  
$\widetilde{H}\subset \widetilde {N}$. It is enough to show that
$\widetilde{G} = \widetilde{H}\cup\gamma_1^{-1}\widetilde{H}$,
i.e., $[\widetilde{G}:\widetilde{H}] = 2$. 

\medskip

\textsl{Step 1.} $\gamma_i^{\pm} \in \gamma_1^{-1}\widetilde{H}$. 

\medskip

We first note that $\gamma_1 = \gamma_1^{-1}\gamma_1^2 \in
\gamma_1^{-1}\widetilde{H}$.  Suppose that $\gamma_i \in \gamma_1^{-1}
\widetilde{H}$ and put $\gamma_i = \gamma_1^{-1}h_i, h_i \in \widetilde{H}$.
Then
\[
\gamma_{i+1} = \gamma_i^{-1}\gamma_i\gamma_{i+1} 
= \gamma_i\gamma_i^{-2}(\gamma_i\gamma_{i+1})
=\gamma_1^{-1}h_i\gamma_i^{-2}(\gamma_i\gamma_{i+1}) \in 
\gamma_i^{-1}\widetilde{H}.
\]
Thus $\gamma_i \in \gamma_1^{-1}\widetilde{H}$ for $i =1, \ldots, 2g+4$.
Also $\gamma_i^{-1} = \gamma_i\gamma_i^{-2} \in \gamma_i^{-1}\widetilde{H}$.

\medskip

\textsl{Step 2.} Any element on $w \in \widetilde{G}$ is in either 
$\widetilde{H}$ or $\gamma_1^{-1}\widetilde{H}$. 

\medskip
Since
\[
\widetilde{G} = \langle \gamma_1,\ldots, \gamma_{2g+4} \mid
\gamma_1\gamma_2\cdot\ldots\cdot\gamma_{2g+4} = 1 \rangle,
\]
we may assume that any $w \in \widetilde{G}$ is of the form
\[
w = w_1\cdot\ldots\cdot w_M,
\]
where $w_k = \gamma_i$ or $\gamma_i^{-1}$ for some $i$.

\medskip

 Suppose that $M$ is even. As $w = (w_1w_2)\cdots(w_iw_{i+1})\cdots(w_{M-1}w_M)$,
 it is enough to show
 \[
 \gamma_i^{\pm}\gamma_j^{\pm}, \gamma_i^{\pm}\gamma_j^{\mp} \in 
 \widetilde{H}
 \]
 for any $i, j$. Since
 \[
 \gamma_i^{-1}\gamma_j^{-1} = (\gamma_j\gamma_i)^{-1}, \quad
 \gamma_i^{-1}\gamma_j = \gamma_i^{-2}\gamma_i\gamma_j, \quad
 \gamma_i\gamma_j^{-1} = \gamma_i\gamma_j\gamma_j^{-2},
 \]
 we only need to show $\gamma_i\gamma_j \in \widetilde{H}$ for any $i, j$.
 There are three possibilities, namely, $(i)$ $i = j$, $(ii)$ $i < j$ and $(iii)$ $i > j$, and
 we check each case separately.
 
 \medskip
 
 The case $(i)$. $\gamma_i\gamma_j = \gamma_i^2 \in \widetilde{H}$.
  
  \medskip
  
  The case $(ii)$. Since
  \[
  \gamma_i\gamma_j = (\gamma_i\gamma_{i+1})\gamma_{i+1}^{-2}(\gamma_{i+1}\gamma_{i+2})\cdot
  \ldots \cdot\gamma_{j-1}^{-2}\gamma_{j-1}\gamma_j,
  \]
  $\gamma_i\gamma_j \in \widetilde{H}$.
  
  \medskip
  The case $(iii)$.  Since
 \[
  \gamma_i\gamma_j = (\gamma_i\gamma_{i-1})\gamma_{i-1}^{-2}(\gamma_{i-1}\gamma_{i-2})\cdot
  \ldots \cdot\gamma_{j+1}^{-2}\gamma_{j+1}\gamma_j,
  \]
  and 
  \[
  \gamma_{k+1}\gamma_k = \gamma_{k+1}^2(\gamma_{k+1}^{-1}\gamma_k^{-1})\gamma_k^2 =
  \gamma_{k+1}^2(\gamma_k\gamma_{k+1})^{-1}\gamma_k^2,
  \]
  we infer that $\gamma_i\gamma_j \in \widetilde{H}$.
  
  Thus $w \in \widetilde{H}$, when $M$ is even.
  
  \medskip
  
  Suppose that $M$ is odd. In this case, $w$ is of the form $\gamma_i\widetilde{w}$ or
  $\gamma_i^{-1}\widetilde{w}$ for some $i$ and $\widetilde{w} \in \widetilde{H}$ by Case 1.
  Since $\gamma_i^{\pm} \in \gamma_1^{-1}\widetilde{H}$, we infer that
  $w \in \gamma_1^{-1}\widetilde{H}$.
  
  From Step 1 and Step 2, we have Claim and Lemma \ref{lem:key-1} follows. \proofend
  
  \subsection{Existence of a $D_6$-cover}
  
  In this section, we show the existence of a normal subgroup $K$ of $G_B$
  such that $(i)$ $K$ is a subgroup of $H$, $(ii)$ $G_B/K \cong D_6$ and
  $(iii)$ the $D_6$-cover corresponding to $K$ is branched along $B$ with
  ramification index $2$.
  
  Let $\mu : \widetilde{S} \to S$ be the canonical resolution of the double
  cover $f : S  \to \PP^1\times\Delta_{\epsilon}$. Since singularities
  of $S$ are on those of $B$, we have $S\setminus f^{-1}(B)
  \cong \tilde{S}\setminus (f\circ\mu)^{-1}(B)$. Hence we have 
  an epimorphism $\overline{\delta} : \pi_1(S\setminus f^{-1}(B), b_+)
  \to \pi_1(\widetilde{S}, b_+)$, where $b_+$ is a point
  in  $f^{-1}(b)$. In particular, we have an epimorphism
  $\delta : H \to H_1(\widetilde{S}, \ZZ)$. 
   Let $A$ be a subgroup of 
  $H_1(\widetilde{S}, \ZZ)$ and put $K_A := \delta^{-1}(A)$. We have
  
  \begin{lem}\label{lem:normal_K}{$K_A$ is a normal subgroup of $G_B$.
  }
  \end{lem}
  
  \proof We first note that $\gamma_i^2 \in \ker\delta$ for $i = 1,\ldots, 2g+4$. Since 
  \[
  \gamma_1^{-1}(\gamma_i\gamma_{i+1})\gamma_1
  = (\gamma_i\gamma_1)^{-1}\gamma_i^2\gamma_{i+1}^2(\gamma_i\gamma_{i+1})^{-1}(\gamma_i\gamma_1),
  \]
  we have
  \[
  \delta(\gamma_1^{-1}(\gamma_i\gamma_{i+1})\gamma_1) = - \delta(\gamma_i\gamma_{i+1}).
  \]
  Let $k$ be an arbitrary element in $K_A$ and suppose that $k$ is of the form
  $u_1\cdot\ldots\cdot u_n$, where $u_i$ ($i = 1,\ldots, n)$ are either $(\gamma_j^2)^{\pm}$
or  $(\gamma_l\gamma_{l+1})^{\pm}$ for some $j, l$. Then
\begin{eqnarray*}
\delta(\gamma_1^{-1}k\gamma_1) & = & \sum_{i=1}^n\delta(\gamma_1^{-1}u_i\gamma_1) \\
    & = & - \sum_{i=1}^n\delta(u_i) \\
    & = & - \delta(k).
 \end{eqnarray*}
 Hence $\gamma_1^{-1}k\gamma_1 \in K_A$. As $H \vartriangleright K_A$ and $G_B = H \cup \gamma_1^{-1}H = 
 H\cup\gamma_1 H$, we have
 $G_B\vartriangleright K_A$. \proofend
 
 \begin{cor}\label{cor:exist_normal} {Let $A$ be a subgroup of $H_1(\widetilde{S}, \ZZ)$ of index $3$.
 Then $K_A$ is a normal subgroup of $G_B$ such that $G_B/K_A \cong D_6$.
 }
 \end{cor}
 
 \proof Chose $h \in H\setminus K_A$ such that $H = K_A\cup hK_A \cup h^{-1}K_A$. Then 
 $\gamma_1^{-1}h\gamma_1 \in h^{-1}K_A$. This implies that $G_B/K_A$ is non-abelian.
 \proofend
 
 \subsection{Proof of Theorem \ref{thm:main}}
 
We keep our notations as before. 
Put $\varphi = \mathrm{pr}_2\circ  f \circ\mu$. 
 We first note that
 $H_1(\widetilde{S}, \ZZ) \cong H_1(|\varphi^{-1}(0)|, \ZZ)$,
 where $|\bullet |$ denotes the underlying topological space of $\bullet$. 
 We call the irreducible component of $\varphi^{-1}(0)$ coming from 
 $\PP^1\times\{0\}$
 \textrm{the main component}. Our first observation is as follows:
 
 \begin{lem}\label{lem:central_fiber}{Let $\overline {\gamma_{2g+2}\gamma_{2g+3}}$ be the class of $\gamma_{2g+2}\gamma_{2g+3}$ in 
 $H_1(\varphi^{-1}(0), \ZZ)$. Then
 \[
 H_1(|\varphi_o^{-1}|(0), \ZZ)\oplus\ZZ\overline{\gamma_{2g+2}\gamma_{2g+3}} \cong H_1(|\varphi^{-1}(0)|, \ZZ).
 \]
 }
 \end{lem}
 
 \proof By observing the difference of the central fibers between $\varphi_o^{-1}(0)$ and
 $\varphi^{-1}(0)$, our statement easily follows. \proofend
 
 \medskip
 
 Put $A = H_1(|\varphi^{-1}(0)|, \ZZ)\oplus \langle
 \overline {(\gamma_{2g+2}\gamma_{2g+3})^3}\rangle$ and let $q : \widehat{X} \to \widetilde{S}$ be the cyclic
 triple cover corresponding to $ K_A$ in $H$, and let $\nu : X \to \hat{S}$ be the
 Stein factorization  of $\widehat{X} \to \widetilde{S} \to \PP^1\times\Delta_{\epsilon}$. 
 By Corollary \ref{cor:exist_normal}, $\hat S$ is a $D_6$-cover and we denote its 
 covering morphism by $\hat {\pi} : \hat S \to \PP^1\times \Delta_{\epsilon}$.
 By our construction,
 $\hat{\pi}^{-1}(\PP^1\times \{0\})$ satisfies the following conditions:
 
 \begin{itemize}
  \item $\hat{\pi}^{-1}(\PP^1\times \{0\})$ consists of three irreducible curves $F_1$, $F_2$ and $F_3$.
  
  \item We may assume that $\sigma^*F_1 = F_1$, $\sigma^*F_2 = F_3$ and
  $\tau^*F_1 = F_2, \tau^*F_2 = F_3$, where $\sigma$ and $\tau$ are
   the elements of $D_6$ as in \S 1.
  
  
  \item For each $i$, $\nu^{-1}(F_i)$ is isomorphic to $\varphi_o^{-1}(0)$.
    The involution induced by $\sigma$ acts on $\nu^{-1}(F_i)$ in the same way as the involution on $\widetilde{W}$ induced by the covering transformation of $f_o : W \to \PP^1\times\Delta_{\epsilon}$ acts on 
    $f_o^{-1}(\PP^1\times\{0\})$.
    
    \end{itemize}
    
   Let $X$ be the quotient surface of $\hat{S}$ by $\sigma$. 
   By Lemma \ref{lem:fund}, $X$ is a non-Galois triple cover of
   $\PP^1\times \Delta_{\epsilon}$ with
   branch locus $B$ and we denote its covering morphism by $\pi :
   X \to \PP^1\times\Delta_{\epsilon}$.  Let $\varphi_X : X \to \Delta_{\epsilon}$ be the induced
   fibration. Then we have
   
  \begin{itemize}
  
  \item $\varphi_X^{-1}(t)$ is a smooth curve of genus $g$ for $t\neq 0$ such that
  $\pi|_{\varphi_X^{-1}(t)} : \varphi_X^{-1}(t) \to \PP^1\times\{t\}$ is a $3$-to-$1$ morphism.
  
  \item $\varphi_X^{-1}(0)$ consists of two reduced components $G_1$ and $G_2$, where $G_1$ is the image of $F_1$ and $G_2$ is the image of
  both $F_2$ and $F_3$. Since $G_1\cong \PP^1$ and $G_1G_2 = 1$,
  $G_1$ is a exceptional curve of the first kind. Also $G_2$ is isomorphic to
  $f_o^{-1}(\PP^1\times\{0\})$. 
  
  \end{itemize}
  Figure 1 explains the case when $g = 3$ and $B_1$ has one $(2, 3)$ cusp.

  Let $x_o$ be a singular point of $X$. $y_o = \pi(x_o)$ is a singular point of $B$ and
  there exists a small neighborhood $U_{y_o}$ and $V_{x_o}$ of $y_o$ and $x_o$, respectively
  such that $\pi(V_{x_o}) = U_{y_o}$ and $\pi|_{V_{x_o}} : V_{x_o} \to U_{y_o}$ is a 
  double cover.
  
  We now blow down $G_1$ and take the canonical resolution $\widetilde{X}$ of all 
  singularities of $X$.
  Since $G_1$ is a vertical divisor, $\varphi_X$ induces another fibration $\varphi_{\widetilde{X}} :
  \widetilde{X} \to \Delta_{\epsilon}$, which gives the desired
  local trigonal fibration in Theorem \ref{thm:main}.

  \begin{figure}[h]
  \centering
\includegraphics[width= 12cm, height=10cm, keepaspectratio]{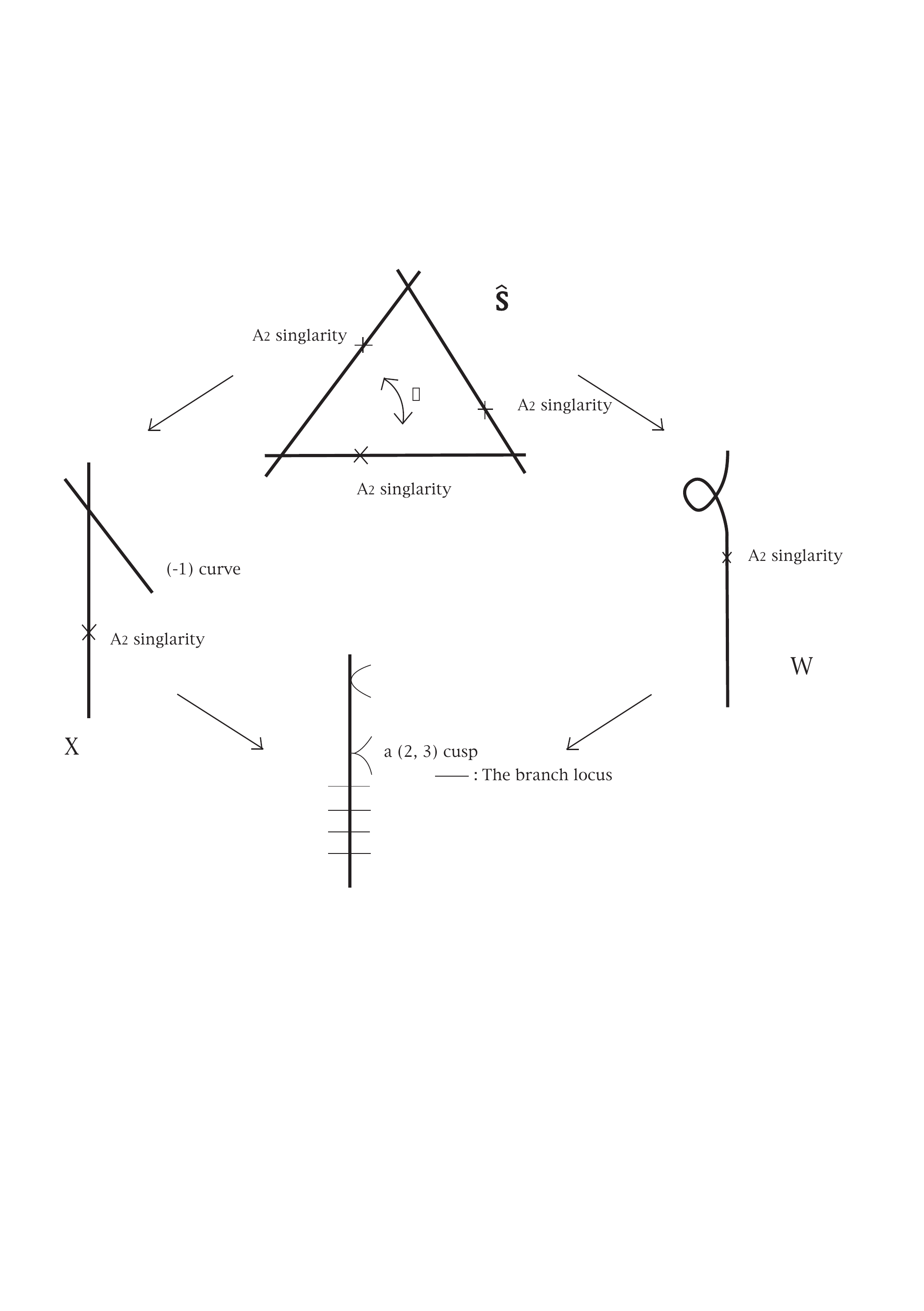}
\caption{The case when  $g = 3$ and  $B_1$ has one $(2, 3)$ cusp.}
\end{figure}

  %
  %


\section{An example for non-horizontal case}

In this section, we give a local trigonal fibration of genus $3$ such that
the central fiber is given as the one for a local hyperelliptic fibration of 
non-horizontal type.

Note that some of local trigonal fibration of non-horizontal type can be reduced to 
those of horizontal type by considering elementary transformations at $\PP^1\times\{0\}$.
Our example is given below is the one which can not be reduced to non-horizontal type.

\begin{exmple}\label{eg:non-horizontal}{\rm
Let $B_1$ be a  reduced divisor on $\PP^1\times \Delta_{\epsilon}$ given by
\[
B_1: t(x^4-t)(x^4+t) = 0,
\]
where $x = X_1/X_0$, $[X_0, X_1]$ being a homogeneous coordinate of $\PP^1$ and $t$ is 
a coordinate of $\Delta_{\epsilon}$. Let $\varphi_o : \caC \to \Delta_{\epsilon}$ be the local
hyperelliptic fibration obtained as in Introduction. $\varphi_o$ is not relatively minimal and
let $\overline{\varphi}_o : \overline{\caC} \to \Delta_{\epsilon}$ be its relatively minimal
model. We denote its central fiber by $F_0$. The configuration of $F_0$ is as follows:
\[
F_0 =4E + 2C_{1} + 2C_{2},
\]
where $E$ is a curve of genus $1$ and $C_{i},  i = 1, 2$ are smooth rational curve with
 $E^2 = -1$, $C_i^2 = -2$,  $EC_{i}=1$ and $C_{1}C_{2}=0$. 
}
\end{exmple}

We first note that $\overline{\varphi}_o : \overline{\caC} \to \Delta_{\epsilon}$ 
in Example~\ref{eg:non-horizontal} is never 
obtained as the relatively minimal model of a local hyperelliptic fibration of horizontal type.
In fact, suppose that there exits a local hyperelliptic fibration of horizontal type such that
the relatively minimal model of 
$\varphi_1 : \caC' \to \Delta_{\epsilon}$ is $\overline{\varphi}_o$.  This means that
$\caC'$ is obtained from $\overline{\caC}$ by a successive blowing-ups. Since $F_0$ has
no reduced component, we infer that the central fiber of $\varphi_1$ has also
no reduced component. On the other hand, the irreducible component of 
the central fiber of $\varphi_1$ arising from $\PP^1\times\{0\}$ is reduced by
taking a local section into account. This leads
us to a contradiction. 

Therefore 
we can not apply Theorem~\ref{thm:main} to obtain a local trigonal fibration with central fiber $F_0$ in Example~\ref{eg:non-horizontal}. 
Nevertheless, there exists a local  trigonal fibration with central fiber $F_0$. We end up this
section in constructing such
an example explicitly.

\begin{exmple}\label{eg:trigonal}{\rm
Let us start with a family of plane quartic curves as follows:

Let $[X_{0},X_{1},X_{2}]$ be a homogeneous coordinates of ${\mathbb P}^{2}$. 
Consider the surface $H_{0}$ of ${\mathbb P}^{2} \times \Delta_{\varepsilon}$ 
defined by 
$$
(X_{1}X_{2}-X_{0}^{2})^{2} + t^{8} (X_{1}^{4}-X_{2}^{4}) = 0.  
$$  

Let  $\pi_{0} \colon H_{0} \rightarrow \Delta_{\varepsilon}$ be the morphism
induced from the second projection 
${\mathbb P}^{2} \times \Delta_{\varepsilon} \rightarrow \Delta_{\varepsilon}$. 
We see that the fiber $\pi_{0}^{-1}(0)$ is the singular locus 
of $H_{0}$. 
Since the fibers $\pi_{0}^{-1}(t)$ $(t \neq 0)$ are nonsingular plane curves of 
degree four, 
they are non-hyperelliptic curves of $3$ three i.e., trigonal curves of genus $3$. 
We blow up ${\mathbb P}^{2} \times \Delta_{\varepsilon}$ along the ideal  
generated by
$X_{1}X_{2} -X_{0}^{2}$ and $t^{4}$. 
Setting $Z_{0}t^{4}=X_{1}X_{2}-X_{0}^{2}$, 
we obtain the defining equation of the proper transformation $H_{1}$ of $H_{0}$ 
in the exceptional set as 
$$
Z_{0}^{2}+(X_{1}^{4}-X_{2}^{4})=0. 
$$
The family $\pi_{1} \colon H_{1} \rightarrow \Delta_{\varepsilon} $ is nonsingular. 
Moreover, since the restriction    
$h \colon \pi_{1}^{-1}(0) \rightarrow \pi_{0}^{-1}(0)$ 
of the morphism $H_{1} \rightarrow H_{0}$ is the double cover  
branched at eight points $X_{1}X_{2}-X_{0}^{2}=X_{1}^{4}-X_{2}^{4}=t=0$, 
the fiber $\pi_{1}^{-1}(0)$ is a hyperelliptic curve.

We now define the automorphism of $G \colon H_{0} \rightarrow H_{0}$ 
induced from the automorphism of ${\mathbb P}^{2} \times \Delta_{\varepsilon}$ as 
$$ 
X_{0} \mapsto iX_{0}, \ \ \ \ X_{1} \mapsto -X_{1}, \ \ \ \ X_{2} \mapsto X_{2}, 
\ \ \ \ t \mapsto it,  
$$ 
where $i=\sqrt{-1}$. 
The fixed points of $G$ are $P_{1}= [0, 0, 1]$ and $P_{2}=[0, 1, 0]$ on 
${\mathbb P}^{2} \times \{ 0 \}$. 
We can naturally define the automorphism $G^{\prime} \colon H_{1} \rightarrow H_{1}$ 
of $H_{1}$ induced from $G$. 
Note that $G^{\prime}$ acts on the coordinate $Z_{0}$ as $Z_{0} \mapsto -Z_{0}$. 
Consider the quotient $S$ of $H_{1}$ by the cyclic group $\left< G^{\prime} \right>$ generated by $G^{\prime}$. 
We see that $(G^{\prime})^4$ is the identity and $(G^{\prime})^{2}$ has four fixed points on 
$\pi_{1}^{-1}(0)$ which are the inverse image of $P_{1}$ and  $P_{2}$ by $h$. 
Since $G^{\prime}$ acts on the coordinate $Z_{0}$ as $Z_{0} \mapsto -Z_{0}$, 
$G^{\prime}$ interchanges the two points of the inverse image of $P_{i}$ $(i=1,2)$ 
by $h$, we see that the restriction $G^{\prime}|_{\pi_{1}^{-1}(0)}$ of $G^{\prime}$ 
is the automorphism of order four and the quotient 
$\pi_{1}^{-1}(0)/ \left< G^{\prime}|_{\pi_{1}^{-1}(0)} \right>$ 
is an elliptic curve. 

Thus, we see that the singular fiber $\phi^{-1}(0)$ of 
$\phi \colon S \rightarrow \Delta_{\varepsilon}$ is the elliptic curve 
with multiplicity four and $S$ has two rational double points of type $A_{1}$ 
on its singular fiber. 
Then, the resolution $\psi \colon \widetilde{S} \rightarrow \Delta_{\varepsilon}$ of the family 
$\phi \colon S \rightarrow \Delta_{\varepsilon}$ has the singular fiber $F_0$ 
and its general fibers are non-hyperelliptic curves. 
}
\end{exmple}

   

Mizuho ISHIZAKA

Graduate School of Mathematical Sciences

University of Tokyo

3-8-1 Komaba Meguro

153-8914 

Tokyo Japan

\bigskip

Hiro-o TOKUNAGA

Department of Mathematics and Information Science

Tokyo Metropolitan University

1-1 Minamiohsawa

 Hachoji 192-0397
 
 Tokyo Japan

\begin{thebibliography}{99}
  %
\bibitem{acc} E.~Artal, J.~Carmona, and J.I.~Cogolludo,
 \emph{Braid monodromy and topology of plane curves}, Duke Math. J.~\textbf{118} 
(2003), no.~2, 261--278.
%

  %
  \bibitem{act} E. ~Artal Bartolo, J.-I. ~Cogolludo and H.~Tokunaga, \emph{
 A survey on Zariski pairs}, to appear in ASPM.
 


  \bibitem{bpv}
W.~Barth, C.~Peters, and A.~Van de~Ven, \emph{Compact complex surfaces}, Erg.
  der Math. und ihrer Grenz., A Series of Modern Surveys in Math.,~\textbf{3}, 
  vol.~4, Springer-Verlag, Berlin, 1984.
  

%
\bibitem{chen-tan} Z.~Chen and S.-L.~Tan, \emph{Upper Bounds on the Slope
of a Genus $3$ Fibration}, Contemp. Math., 400, Amer. Math.Soc.



%
\bibitem{dimca}
A.~Dimca, \emph{Singularities and topology of hypersurfaces}, Springer-Verlag,
  New York, 1992.
%

\bibitem{sga}
A.~Grothendieck, \emph{Rev\^etements \'etales et groupe fondamental}, Lecture
  Notes in Math., \textbf{224}, Springer-Verlag, Berlin, 1971.

  \bibitem{horikawa} 
E.~Horikawa, \emph{On deformation of quintic surfaces},
Invent. Math.~\textbf{31} (1975), \rm 43--85.

%


\bibitem{shimada-tokunaga} I.~Shimada and H.~Tokunaga, \emph{
The fundamental group and singluarities (in japanese)}, 
Dais\^ukyokusen to tokuiten, Kyoritsu Shuppan.
%

\bibitem{tokunaga2} H.~Tokunaga,
\emph{Triple coverings of algebraic surfaces accoring to the Cardano formula}, 
J. of Math. Kyoto Univ.~\textbf{31} (1991), 359--375.

  \end{thebibliography}
\end{document}